\documentclass{article}
\usepackage[utf8]{inputenc}
\usepackage[margin=1in, includehead, includefoot]{geometry}
\usepackage{amsthm,amssymb,amsmath}
\usepackage{enumitem}
\usepackage{mathrsfs}
\usepackage{mathtools}
\usepackage{fancyhdr}
\usepackage{indentfirst}
\usepackage{graphicx}
\usepackage{placeins}
\usepackage{wrapfig}
\usepackage{url}
\usepackage[numbers]{natbib}
\usepackage[hidelinks]{hyperref}

\usepackage{xcolor}

\newtheorem{definition}{Definition}
\newtheorem{lemma}{Lemma}[section]
\newtheorem{theorem}{Theorem}[section]
\newtheorem{corollary}{Corollary}[section]
\newtheorem{observation}{Observation}[section]
\newtheorem{proposition}{Proposition}[section]

\newcommand{\toroman}[1]{\textit{\expandafter{\romannumeral #1\relax}}}

\newcommand{\cbeginproof}[0]{\par\noindent\textit{Proof.} }
\newcommand{\cendproof}[0]{ \qed\par\vspace{1em}}

\newcommand{\npcompleteproblem}[3]{ \par\vspace{0.5em}\noindent{\textbf{#1}}\newline \textbf{INSTANCE: } #2 \newline \textbf{QUESTION: } #3 \par\vspace{0.5em} }

\title{Open-locating-dominating sets with error correction}
\author{
    \small Devin C. Jean\\
    \small Computer Science Department \\
    \small Vanderbilt University\\
    \small \texttt{devin.c.jean@vanderbilt.edu}
    \and
    \small Suk J. Seo\\
    \small Computer Science Department\\
    \small Middle Tennessee State University\\
    \small \texttt{Suk.Seo@mtsu.edu}
}
\date{}

\begin{document}
\maketitle
\thispagestyle{empty}

\begin{abstract}
An open-locating-dominating set of a graph models a detection system for a facility with a possible ``intruder" or a multiprocessor network with a possible malfunctioning processor.
A ``sensor'' or ``detector'' is assumed to be installed at a subset of vertices where it can detect an intruder or a malfunctioning processor in their neighborhood, but not at itself.
We consider a fault-tolerant variant of an open-locating-dominating set called an error-correcting open-locating-dominating set, which can correct a false-positive or a false-negative signal from a detector.
In particular, we prove the problem of finding a minimum error-correcting open-locating-dominating set in an arbitrary graph is NP-complete.
Additionally, we characterize the existence criteria for an error-correcting open-locating-dominating sets for an arbitrary graph.
We also consider extremal graphs that require every vertex to be a detector and minimum error-correcting open-locating-dominating sets in infinite grids.
\end{abstract}


\section{Introduction}\label{sec:intro}

Let $G$ be a graph with vertices $V(G)$ and edges $E(G)$.
The \emph{open-neighborhood} of a vertex $v \in V(G)$, denoted $N(v)$, is the set of all vertices adjacent to $v$: $\{w \in V(G) : vw \in E(G)\}$.
The \emph{closed-neighborhood} of a vertex $v \in V(G)$, denoted $N[v]$, is the set of all vertices adjacent to $v$, as well as $v$ itself: $N(v) \cup \{v\}$.
A vertex set $S \subseteq V(G)$ is a dominating set if all vertices are within the closed neighborhood of some $v \in S$; that is, $\cup_{v \in S}{N[v]} = V(G)$.
Similarly, $S \subseteq V(G)$ is an open-dominating set if all vertices are within the open neighborhood of some $v \in S$; that is, $\cup_{v \in S}{N(v)} = V(G)$.

Assume that $G$ models a facility with a possible ``intruder'' (or a multiprocessor network with a possible malfunctioning processor), where we want to be able to precisely determine the location of the intruder by placing (the minimum number of) detectors. Much work has been done on the topic of the graphical parameters involving various types of detectors.
If a detector at $v$ can determine whether if the intruder is in $N[v]$, then the detection system is called an identifying code (IC) as introduced by Karpovsky, Charkrabarty and Levitin \cite{karpovsky}.
A subset $S$ of $V(G)$ is an \textit{identifying code} if it is a dominating set for which, given any two vertices $u$ and $v$ in $V(G)$, we have $N[u] \cap S \neq N[v] \cap S$. 
If a detector at $v$ can determine whether if the intruder is in $N[v]$ and can tell whether the intruder is at $v$ or in $N(v)$, then the detection system is called locating-dominating (LD) set, as introduced in Slater \cite{dom-loc-acyclic}. 
A subset $S$ of $V(G)$ is a \textit{locating-dominating set} if it is a dominating set for which, given any two vertices $u$ and $v$ in $V(G) - S$, we have $N(u) \cap S \neq N(v) \cap S$.
Lobstein \cite{dombib} maintains a bibliography, currently with more than 500 entries, for work on detection system related topics.

This paper focuses on a detection system where a sensor at vertex $v$ can determine if an intruder is in $N(v)$, but notably not at $v$ itself.
This type of detection system is called an open-locating-dominating (OLD) set, as introduced by Honkala, Laihonen and Ranto \cite{honk02d} and Seo and Slater \cite {old}.
A subset $S$ of $V(G)$ is an \textit{open-locating-dominating set} if it is an open-dominating set for which, given any two vertices $u$ and $v$ in $V(G)$, we have $N(u) \cap S \neq N(v) \cap S$.

For an OLD set $S \subseteq V(G)$ and $u \in V(G)$, we let $N_S(u) = N(u) \cap S$ denote the (open) \emph{dominators} of $u$ and $dom(u) = |N_S(u)|$ denote the (open) \emph{domination number} of $u$.
A vertex $v \in V(G)$ is \emph{$k$-open-dominated} by an open-dominating set $S$ if $|N_S(v)| = k$.
If $S$ is an open-dominating set and $u,v \in V(G)$, $u$ and $v$ are \emph{$k$-distinguished} if $|N_S(u) \triangle N_S(v)| \ge k$, where $\triangle$ denotes the symmetric difference.
If $S$ is an open-dominating set and $u,v \in V(G)$, $u$ and $v$ are \emph{$k^\#$-distinguished} if $|N_S(u) - N_S(v)| \ge k$ or $|N_S(v) - N_S(u)| \ge k$.
We will also use terms such as ``at least $k$-dominated'' to denote $j$-dominated for some $j \ge k$.

There are several fault-tolerant variants of OLD sets.
For example, a redundant open-locating-dominating set is resilient to a detector being destroyed or going offline \cite{ftold}.
Thus, an open-dominating set $S \subseteq V(G)$ is called a \emph{redundant open-locating-dominating (RED:OLD)} set if $\forall v \in S$, $S - \{v\}$ is an OLD set.
Another variant of an OLD set  called an \textit{error-detecting open-locating-dominating} (DET:OLD) set, which is capable of correctly identifying an intruder even when at most one sensor or detector incorrectly reports that there is no intruder \cite{ftold}. 
Hence, DET:OLD sets allow for uniquely locating an intruder in a way which is resilient to up to one false negative.

The focus of this paper is another variant of an OLD set  called an \textit{error-correcting open-locating-dominating} (ERR:OLD) set, which is capable of correcting false negative or false positive errors.
The following theorem characterizes OLD, RED:OLD, DET:OLD, and ERR:OLD sets and they are useful in constructing those sets or verifying whether a given set meets their requirements.

\begin{theorem}\label{theo:multi-param} 
An open-dominating set is
\begin{enumerate}[noitemsep, label=\roman*.]
    \item \cite{ourtri} an OLD set if and only if every pair of vertices is 1-distinguished.
    \item \cite{ftold} a RED:OLD set if and only if all vertices are at least 2-dominated and all pairs are 2-distinguished.
    \item \cite{ftold} a DET:OLD set if and only if all vertices are at least 2-dominated and all pairs are $2^\#{}$-distinguished.
    \item \cite{ftold} an ERR:OLD set if and only if all vertices are at least 3-dominated and all pairs are 3-distinguished.
\end{enumerate}
\label{theo:old-redold-detold}
\end{theorem}

\begin{figure}[ht]
    \centering
    \begin{tabular}{c@{\hspace{2.5em}}c@{\hspace{2.5em}}c@{\hspace{2.5em}}c}
        \includegraphics[width=0.2\textwidth]{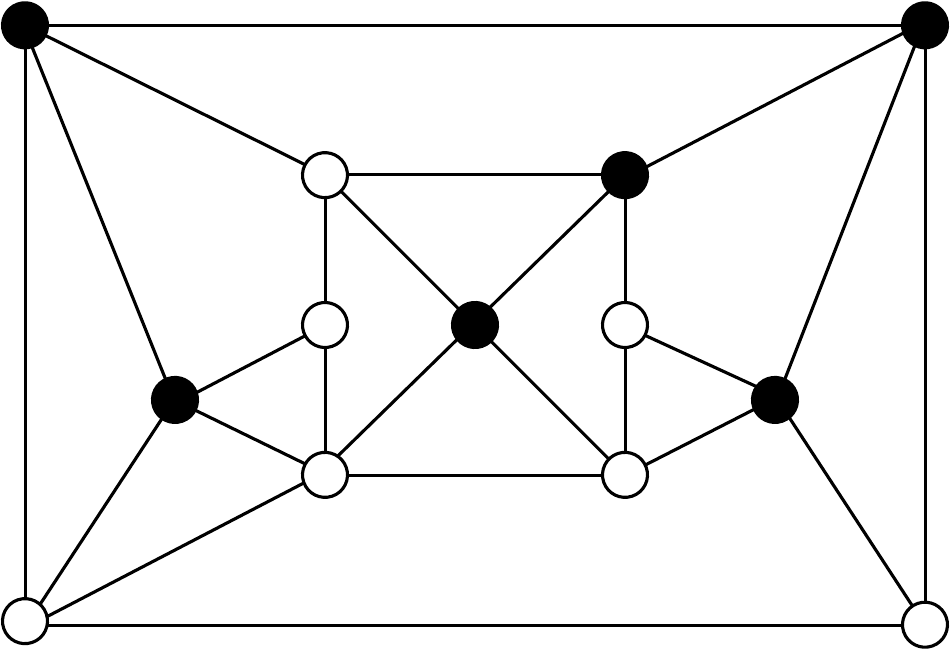} & \includegraphics[width=0.2\textwidth]{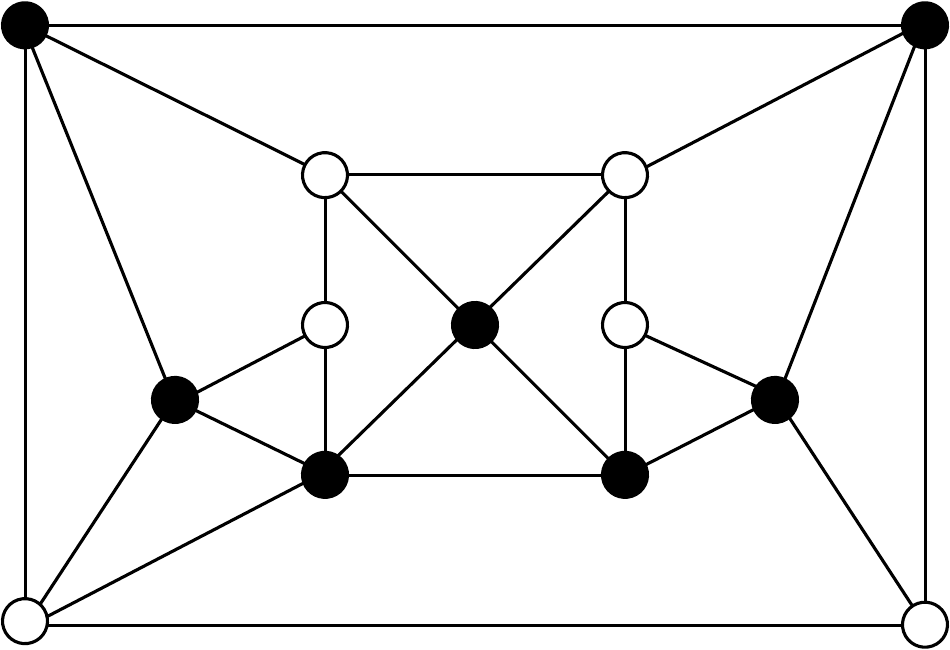} & \includegraphics[width=0.2\textwidth]{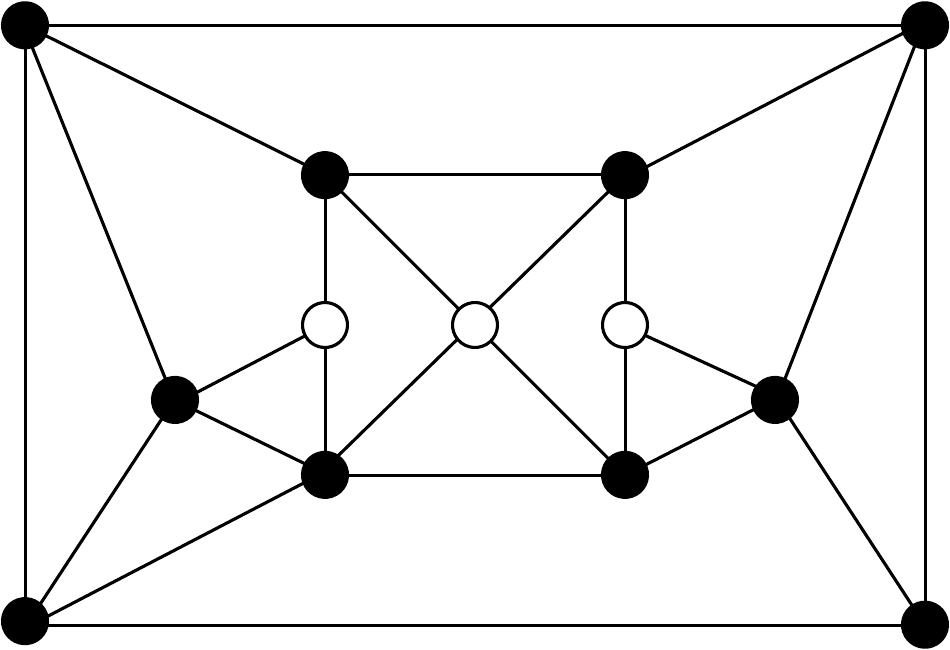} & \includegraphics[width=0.2\textwidth]{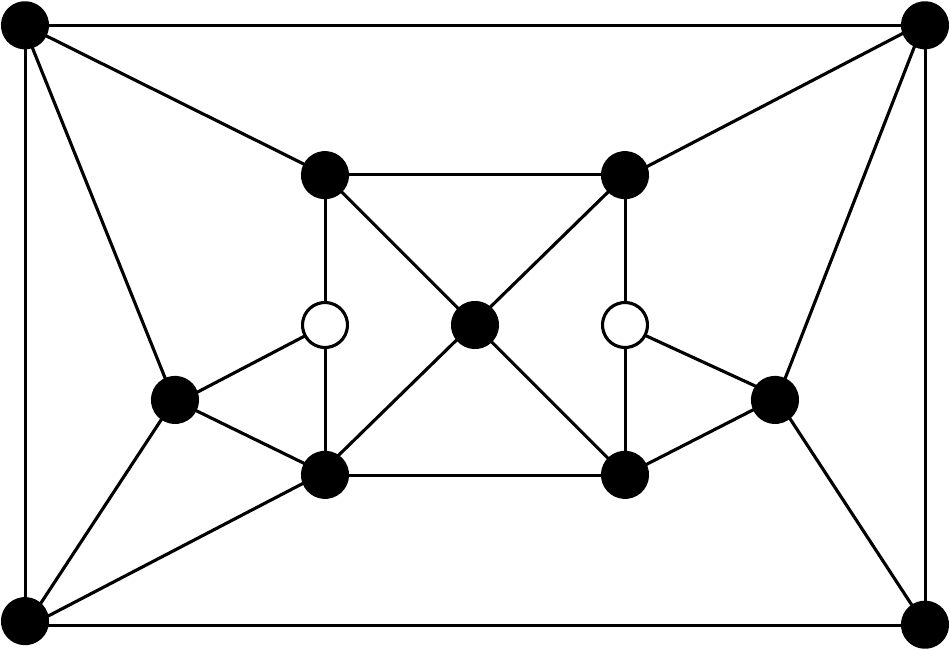} \\
        (a) & (b) & (c) & (d)
    \end{tabular}
    \caption{Optimal OLD (a), RED:OLD (b), DET:OLD (c), and ERR:OLD (d) sets for $G_{13}$. Darkened vertices represent detectors.}
    \label{fig:ex-finite}
\end{figure}

For finite graphs, the notations $\textrm{OLD}(G)$, $\textrm{RED:OLD}(G)$, $\textrm{DET:OLD}(G)$, and $\textrm{ERR:OLD}(G)$ represent the cardinality of the smallest possible OLD, RED:OLD, DET:OLD, and ERR:OLD sets on graph $G$, respectively \cite{ftold}.
From Figure~\ref{fig:ex-finite}, which shows optimal solutions on the given graph which we will call $G$, we see that $\textrm{OLD}(G) = 6$, $\textrm{RED:OLD}(G) = 7$, $\textrm{DET:OLD}(G) = 10$, and $\textrm{ERR:OLD}(G) = 11$.

In Section~\ref{sec:errold-char}, we characterize the existence criteria for an error-correcting open-locating-dominating sets for an arbitrary graph.
We also consider extremal graphs that require every vertex to be a detector.
In Section~\ref{sec:npc}, we prove the problem of finding a minimum error-correcting open-locating-dominating set in an arbitrary graph is NP-complete.
And finally, in Section~\ref{sec:inf-grid}, we investigate on minimum error-correcting open-locating-dominating sets in infinite grids.

\section{Existence and extremal graphs for ERR:OLD}\label{sec:errold-char}

Because an ERR:OLD set requires 3-domination, many graphs do not allow an ERR:OLD set to exist.
In light of this, our next theorem characterizes all graphs that permit an ERR:OLD set.

\begin{theorem}\label{theo:err-old-exist}
A graph, $G$, permits an ERR:OLD set if and only if $\delta(G) \ge 3$ and for every $C_4$ subgraph, $abcd$, $|N(a) \triangle N(c)| \ge 3$.
\end{theorem}
\begin{proof}
First, assume to the contrary that there exists a vertex $v$ with $deg(v) \le 2$; then, $v$ cannot be 3-dominated, a contradiction.
Similarly, if we assume there exists $a,c \in V(G)$ with $|N(a) \triangle N(c)| \le 2$, then $a$ and $c$ could not be distinguished, a contradiction.

For the converse, suppose $G$ satisfies both properties of the theorem statement.
To demonstrate the existence of ERR:OLD, we will conservatively let $S = V(G)$ be the detector set.
Clearly, $\delta(G) \ge 3$ guarantees all vertices will be 3-dominated, so we need only show that all vertex pairs are distinguished.
Let $u,v \in V(G)$ with $u \neq v$.
If $d(u,v) \ge 3$ then $u$ and $v$ will have non-intersecting neighborhoods and thus be distinguished; otherwise, we can assume $d(u, v) \le 2$.
Suppose $d(u, v) = 2$, and let $uxv$ be a path from $u$ to $v$.
If $|N(u) \cap N(v)| = 1$, then $u$ and $v$ must be at least 4-distinguished, and we would be done, so assume $\{x,y\} \subseteq N(u) \cap N(v)$ with $x \neq y$.
Then we have a $C_4$ subgraph, $uxvy$, so by assumption $|N(u) \triangle N(v)| \ge 3$ and $(u, v)$ are distinguished.
Finally, suppose $d(u, v) = 1$.
If $(N(u) \triangle N(v)) - \{u,v\} \neq \varnothing$, then $u$ and $v$ will be 3-distinguished, and we would be done; otherwise, we assume that $N(u) \triangle N(v) = \{u,v\}$.
Because $\delta(G) \ge 3$, it must be that $\exists p,q \in N(u) \cap N(v)$ with $p \neq q$.
Then we have a $C_4$ subgraph $upvq$, so by assumption $u$ and $v$ must be distinguished, completing the proof.
\end{proof}

\begin{corollary}
A $C_4$-free graph, $G$, permits an ERR:OLD set if and only if $\delta(G) \ge 3$.
\end{corollary}

Two distinct vertices $u,v \in V(G)$ are said to be \emph{twins} if $N[u] = N[v]$ (\emph{closed twins}) or $N(u) = N(v)$ (\emph{open twins}).
From Theorem~\ref{theo:err-old-exist}, we see that if $G$ has twins $u$ and $v$, then it cannot have an ERR:OLD set since $deg(u) = deg(v) \le 2$ or $u$ and $v$ form the corners of a $C_4$ subgraph with $|N(u) \triangle N(v)| \le 2$.

\begin{corollary}
If $G$ permits an ERR:OLD set, then $G$ is twin-free.
\end{corollary}

Next, we will show that there are no graphs with ERR:OLD sets on $n \le 6$ vertices.

\begin{theorem}\label{theo:err-old-n7}
If $G$ has an ERR:OLD, then $n \ge 7$.
\end{theorem}
\begin{proof}
Firstly, we know that $G$ cannot be acyclic because an ERR:OLD set requires $\delta(G) \ge 3$.
Thus, $G$ has a cycle, which we will exploit in casing.
We will proceed by casing on the size of the smallest cycle, which implies that edges cannot be added between the vertices of said smallest cycle.
We assume $n \le 6$, as otherwise we would be done.
If the smallest cycle is $C_k$ for $5 \le k \le 6$, then every vertex in the cycle must have a third dominator, and these additional $k$ vertices must all be distinct to avoid creating smaller cycles; thus $n \ge 2k \ge 10 > 6$, a contradiction.
Suppose the smallest cycle is a $C_4$ subgraph $abcd$, then each vertex in the cycle must have a third dominator, $a'$, $b'$, $c'$, and $d'$ adjacent to $a$, $b$, $c$, and $d$, respectively.
We see that $a' \neq b' \neq c' \neq d' \neq a'$; however, because $n \le 6$ it must be that $a' = c'$ and $b' = d'$.
We see that $a$ and $c$ cannot have any additional edges without creating a smaller cycle, but they are open twins and so cannot be distinguished, a contradiction.
Otherwise we can assume that $G$ has a triangle.

Suppose we have $\{b,d\} \subseteq N(a) \cap N(c)$ and $bd \in E(G)$, which we call a ``diamond" subgraph.
If $ac \notin E(G)$ then distinguishing $a$ and $c$ will result in $n \ge 7$, a contradiction; thus, we can assume that $ac \in E(G)$, meaning we have a $K_4$ subgraph.
By symmetry, we can assume that $a$, $b$, and $c$ have third dominators $a'$, $b'$, and $c'$.
If $a'$, $b'$, and $c'$ are all distinct, then $n \ge 7$ and we would be done; otherwise without loss of generality assume $a' = b'$.
Suppose $a' = b' = c'$; then we require $da' \in E(G)$ and by symmetry $\exists x \in N(d) - N(a') - \{a,b,c,d,a'\}$.
We now have a $K_5$ subgraph plus an extra vertex $x \in N(d)$.
Distinguishing $a$ and $b$ requires without loss of generality $x \in N(a)$, but then $a$ and $d$ cannot be distinguished, a contradiction.
Otherwise, we can assume $a' = b' \neq c'$; we see that distinguishing $a$ and $b$ requires, without loss of generality, $bc' \in E(G)$.
If $a'd \in E(G)$ then $dc' \in E(G)$ to distinguish $a$ and $d$; however, $b$ and $d$ can no longer be distinguished, a contradiction, so we assume $a'd \notin E(G)$ and $dc' \notin E(G)$.
If $a'c \in E(G)$ then $b$ and $c$ cannot be distinguished, a contradiction, so by symmetry we assume $a'c \notin E(G)$ and $c'a \notin E(G)$.
Thus, 4-dominating $a'$ and $c'$ requires $a'c' \in E(G)$, but then $c$ and $a'$ cannot be distinguished, a contradiction.
Thus, the existence of a diamond subgraph implies that $n \ge 7$.

At this point, we know that $G$ has a triangle, say $abc$.
To 3-dominate the vertices in $abc$, assume $\exists a' \in N(a)-\{b,c\}$, $\exists b' \in N(b)=\{a,c\}$, and $c' \in N(c)-\{a,b\}$.
Additionally, we can assume $a'$, $b'$, and $c'$ are distinct because otherwise we produce a diamond subgraph.
To 3-dominate $a'$, $b'$, and $c'$ without producing a diamond subgraph, it must be that $\{a'b', b'c', c'a'\} \subseteq E(G)$.
At this point we can no longer add any edges without producing a diamond subgraph, but $a'$ and $b$ are not distinguished, a contradiction, completing the proof.
\end{proof}

\begin{figure}[ht]
    \centering
    \begin{tabular}{c@{\hspace{6em}}c}
        \includegraphics[width=0.25\textwidth]{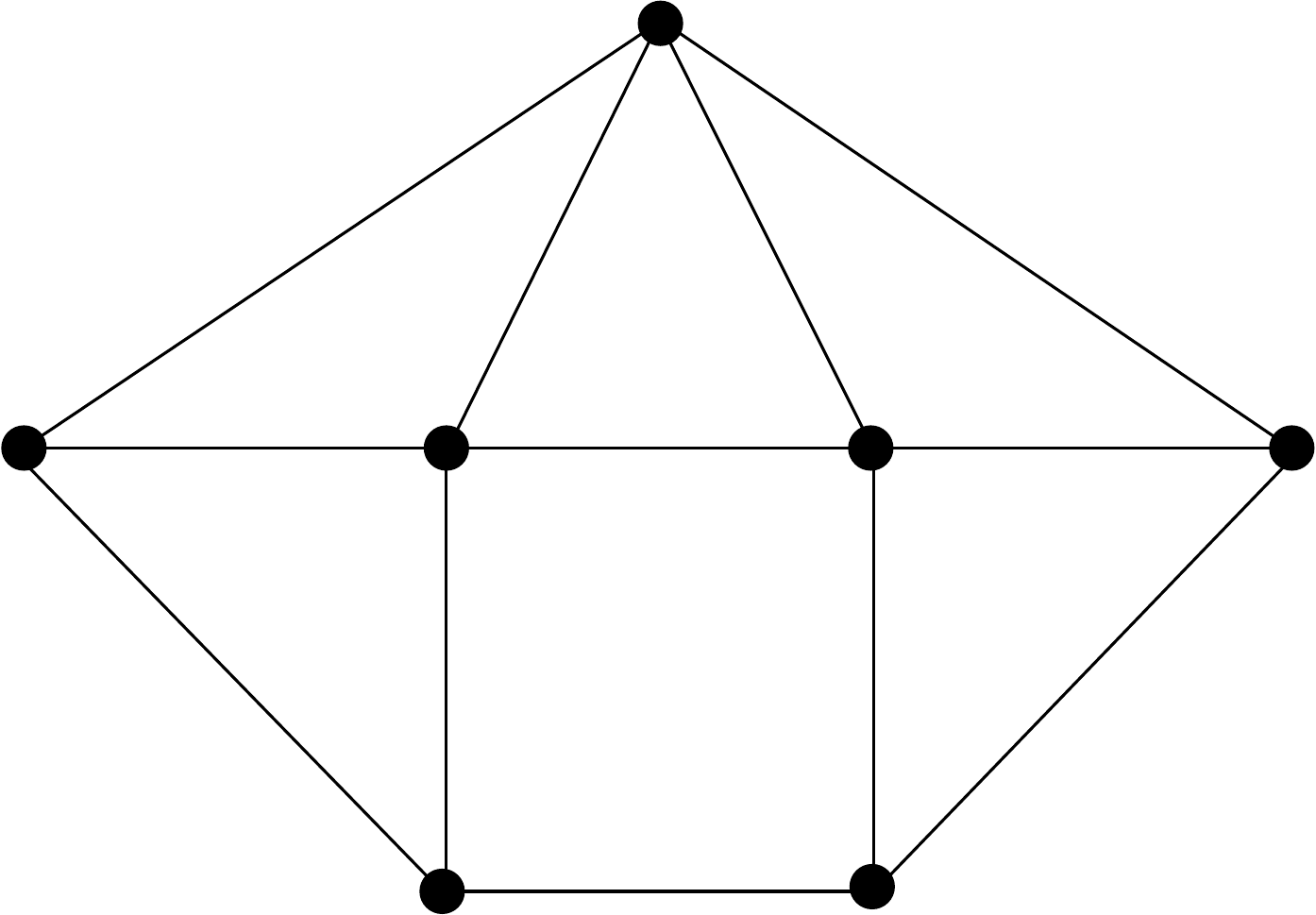} & \includegraphics[width=0.2\textwidth]{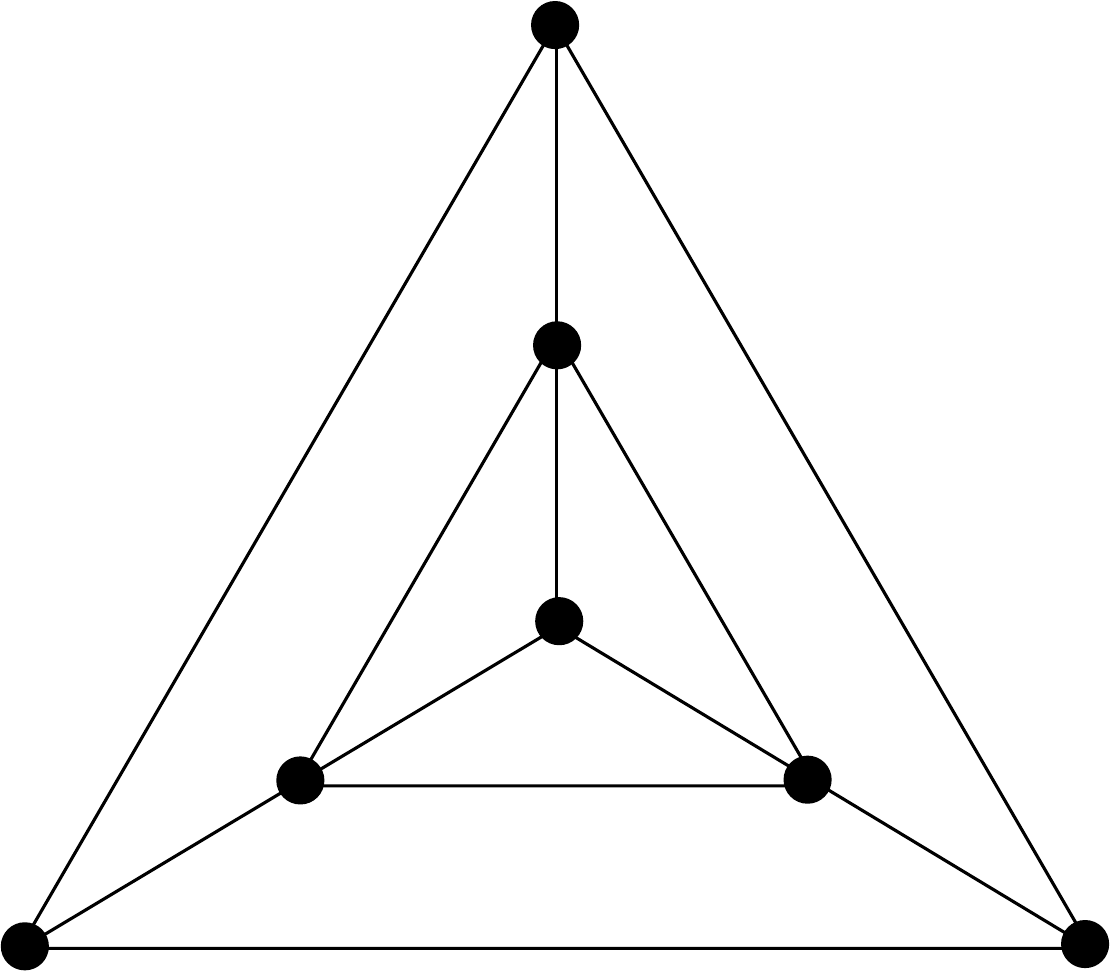} \\
        (a) & (b)
    \end{tabular}
    \caption{The two smallest graphs supporting ERR:OLD, which have $n = 7$ and $m = 12$}
    \label{fig:g7}
\end{figure}

Figure~\ref{fig:g7} shows the first graphs with an ERR:OLD set in the lexicographic ordering of $(n,m)$ tuples; i.e., the graphs with the smallest number of edges given the smallest number of vertices.
For ERR:OLD, this is $n = 7$ and $m = 12$.

Next, we will show the relation between the existence of ERR:OLD sets and DET:OLD sets for regular graphs.

\begin{observation}\label{prop:k-sharp-2k}
If $A$ and $B$ are sets with $|A| = |B|$, then $|A - B| = k$ implies $|A \triangle B| = 2k$.
\end{observation}
\begin{proof}
Suppose $|A - B| = k$.
Then $|A \cap B| = |A| - k = |B| - k$.
Therefore $|A \triangle B| = |A - B| + |B - A| = k + (|B| - |B \cap A|) = k + (|B| - (|B| - k)) = k + k$, completing the proof.
\end{proof}

Suppose $G$ is a regular graph with DET:OLD.
We can conservatively let $S = V(G)$ be the set of detectors.
Then every vertex pair in $G$ must be $2^{\#}$-distinguished; thus, Proposition~\ref{prop:k-sharp-2k} yields that every vertex pair is $4$-distinguished, which is sufficient for the distinguishing requirement of ERR:OLD.

\begin{proposition}
A regular graph $G$ with $\delta(G) = \Delta(G) \ge 3$ has DET:OLD if and only if it has ERR:OLD.
\end{proposition}

A graph is cubic or 3-regular if every vertex is of degree 3.
Jean and Seo \cite{our-18} have previously proven that a cubic graph has DET:OLD if and only if the graph is $C_4$-free; thus, we have the following existence criteria.

\begin{corollary}\label{cor:c4-free}
A cubic graph $G$ has an ERR:OLD set if and only if $G$ is $C_4$-free.
\end{corollary}

\begin{figure}[ht]
    \centering
    \begin{tabular}{c@{\hspace{6em}}c}
        \includegraphics[width=0.3\textwidth]{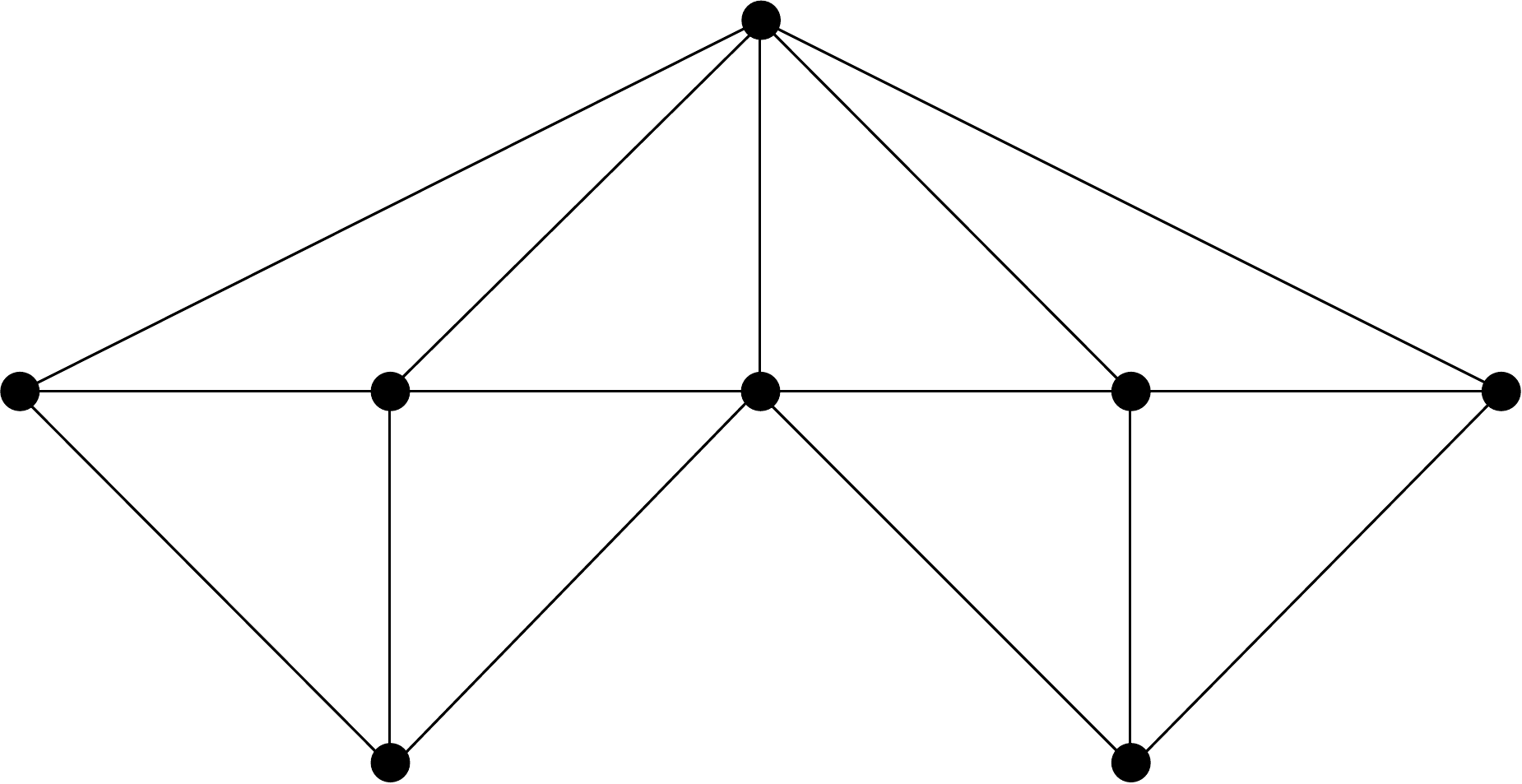} & \includegraphics[width=0.2\textwidth]{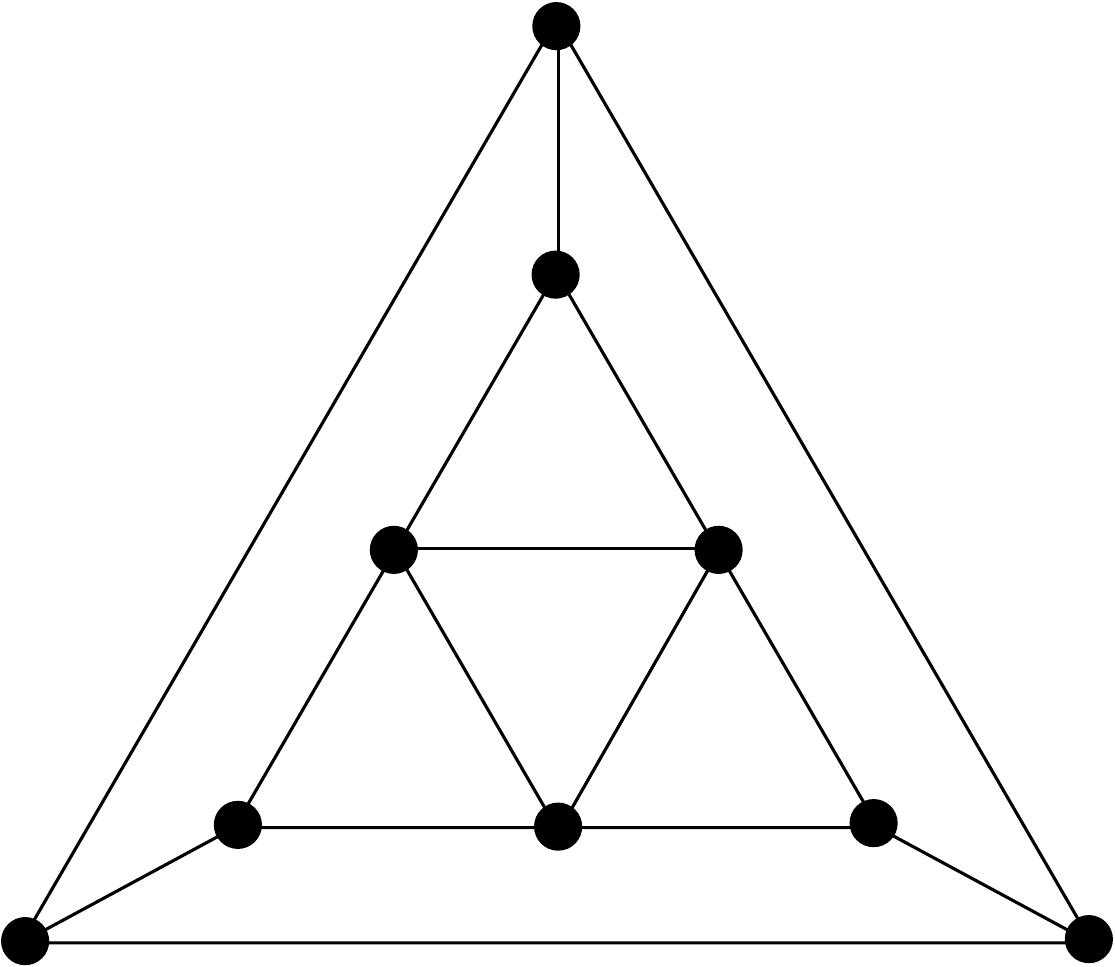} \\
        (a) & (b)
    \end{tabular}
    \caption{Extremal graphs with $\textrm{ERR:OLD}(G) = n$ with the fewest number of edges when $n = 8$ and $n = 9$}
    \label{fig:g8-9}
\end{figure}

Now we consider extremal graphs with $\textrm{ERR:OLD}(G) = n$ with least $m$.
If a graph $G$ has an ERR:OLD set, then Theorem~\ref{theo:err-old-exist} yields that $\delta(G) \ge 3$, which serves as a lower bound for the least $m$, namely $m \ge \frac{3}{2}n$.
Figures \ref{fig:g7} and \ref{fig:g8-9} show extremal graphs with $\textrm{ERR:OLD}(G) = n$ with the fewest number of edges when $7 \le n \le 9$.
For $n \ge 10$, we will show that cubic graphs exist which meet the lower bound of $\delta(G) = \Delta(G) = 3$; however, cubic graphs do not exist when $n$ is odd, so we introduce a new class of almost-cubic graphs:

\begin{definition}
$G$ is \emph{quasi-cubic} if every vertex is degree 3 aside from one degree-4 vertex.
\end{definition}

When $n \ge 10$, cubic and quasi-cubic graphs together form the entire family of extremal graphs with $\textrm{ERR:OLD}(G) = n$ which have the minimum number of edges for a given $n$.
If $S$ is an ERR:OLD set, we know that $v \in S$ if $\exists w \in N(v)$ with $deg(w) = 3$.
Because every vertex in a cubic or quasi-cubic graph is adjacent to a degree-3 vertex, these classes or graphs require $\textrm{RED:OLD}(G) = n$ if they permit ERR:OLD.

\begin{theorem}\label{theo:errold-cubic-quasi-cubic-c4-exist}
If $G$ is cubic or quasi-cubic, then $G$ has an ERR:OLD set if and only if $G$ is $C_4$-free.
\end{theorem}
\begin{proof}
Suppose $G$ has a $C_4$ subgraph, $abcd$.
Because $G$ is cubic or quasi-cubic, at most one vertex in $abcd$ can be degree 4, and all others will be degree 3.
By symmetry, assume $deg(a)=deg(c)=3$.
Then $a$ and $c$ cannot be 3-distinguished, implying $G$ does not have ERR:OLD.

For the converse, suppose $G$ does not have ERR:OLD.
Because $\delta(G) \ge 3$, we know that an ERR:OLD set not existing implies $\exists u,v \in V(G)$ where $u$ and $v$ are not 3-distinguished.
Regardless of whether or not $uv \in E(G)$, we see the only way to have $u$ and $v$ not be 3-distinguished is by having $|N(u) \cap N(v)| \ge 2$, which creates a $C_4$ subgraph, completing the proof.
\end{proof}

Suppose we have a cubic graph $G$ that permits an ERR:OLD set.
Then, we can construct a quasi-cubic graph from $G$ as shown in the next theorem. 
 
\begin{theorem}
Let $G$ be a cubic graph with ERR:OLD where $a,b,c,d \in V(G)$ and $ab,cd \in E(G)$ such that $ab$ and $cd$ are not part of a triangle or the terminal edges of a $P_5$ subgraph.
Construct quasi-cubic $G'$ from $G$ by deleting edges $ab$ and $cd$ and adding a new vertex $x$ with $N(x) = \{a,b,c,d\}$.
Then $G'$ has ERR:OLD.
\end{theorem}
\begin{proof}
By Theorem~\ref{theo:errold-cubic-quasi-cubic-c4-exist}, we know that $G$ is $C_4$-free and we need only show that $G'$ is likewise $C_4$-free.
Suppose for a contradiction that $G'$ has a $C_4$ subgraph; then due to $G$ being $C_4$-free, the $C_4$ subgraph must involve vertex $x$, which is an endpoint of every new edge in $G'$.
By symmetry, we can reduce the problem to two possible cases: the $C_4$ includes edges $\{ax,xb\}$ or edges $\{ax,xc\}$, and we will consider $y \in V(G)$ to be the opposite vertex of $x$ in the $C_4$ subgraph.
If we are in the $\{ax,xb\}$ case, we see that any choice of $y$ results in $aby$ forming a triangle in $G$, a contradiction.
Otherwise, in the $\{ax,xc\}$ case we see that $y \neq b$ and $y \neq d$ because $ab,cd \notin E(G')$, so $y$ causes $ab$ and $cd$ to be the ends of a $P_5$ subgraph in $G$, a contradiction, completing the proof.
\end{proof}

\section{NP-completeness of Error-correcting OLD}\label{sec:npc}

It has been shown that many graphical parameters related to detection systems, such as finding the cardinality of the smallest IC, LD, or OLD sets, are NP-complete problems \cite{ld-ic-np-complete-2, NP-complete-ic, NP-complete-ld, old}.
Similarly, the problems involving RED:OLD and DET:OLD parameters have been known to be NP-complete \cite{redold, our-18}.
We will now prove that the problem of determining the smallest ERR:OLD set is also NP-complete.
For additional information about NP-completeness, see Garey and Johnson \cite{np-complete-bible}.

\npcompleteproblem{3-SAT}{Let $X$ be a set of $N$ variables.
Let $\psi$ be a conjunction of $M$ clauses, where each clause is a disjunction of three literals from distinct variables of $X$.}{Is there is an assignment of values to $X$ such that $\psi$ is true?}

\npcompleteproblem{Error-Detecting Open-locating-dominating set (ERR-OLD)}{A graph $G$ and integer $K$ with $7 \le K \le |V(G)|$.}{Is there an ERR:OLD set $S$ with $|S| \le K$? Or equivalently, is ERR:OLD($G$) $\le K$?}

\begin{figure}[ht]
    \centering
    \includegraphics[width=0.2
    \textwidth]{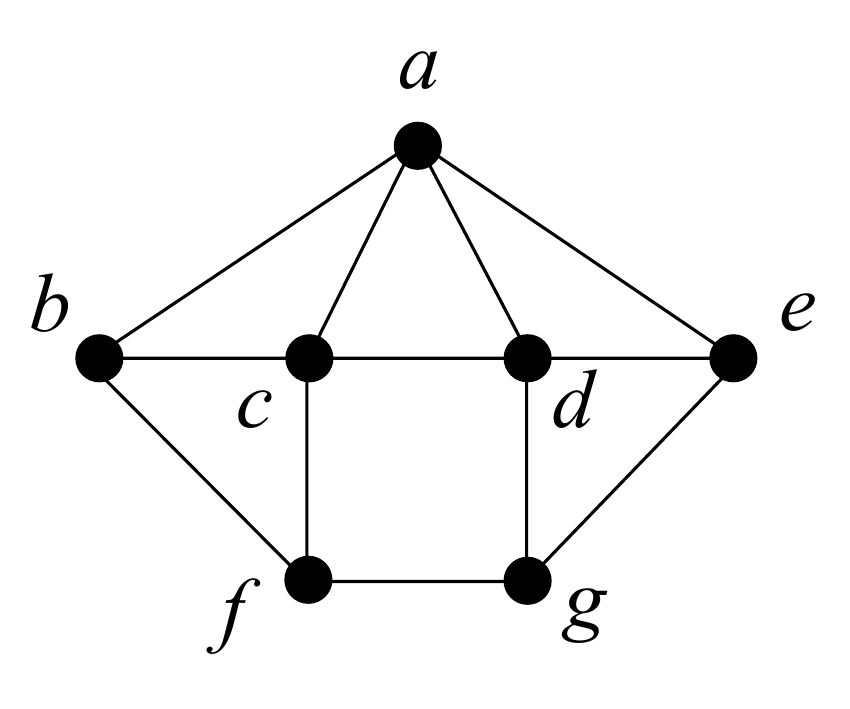}
    \caption{Subgraph $G_7$}
    \label{fig:g7-labeled}
\end{figure}

\begin{lemma}\label{lem:g7}
Given a graph $G$, let $G_{7}$ be the subgraph shown in Figure~\ref{fig:g7-labeled}, where vertices $\{b,c,d,e\}$ have no additional edges in the larger graph but $\{a,f,g\}$ are permitted to have additional edges to vertices outside of this subgraph.
If $S$ is an ERR:OLD set for $G$, then $V(G_{7}) \subseteq S$.
\end{lemma}
\begin{proof}
To 3-distinguish $b$ and $c$, we require $\{b,c,d\} \subseteq S$.
To 3-dominate $b$, we require $\{a,f\} \subseteq S$.
By symmetry, we have that $V(G_{7}) \subseteq S$, completing the proof.
\end{proof}

\begin{theorem}
The ERR-OLD problem is NP-complete.
\end{theorem}
\begin{wrapfigure}{r}{0.4\textwidth}
    \centering
    \includegraphics[width=0.4\textwidth]{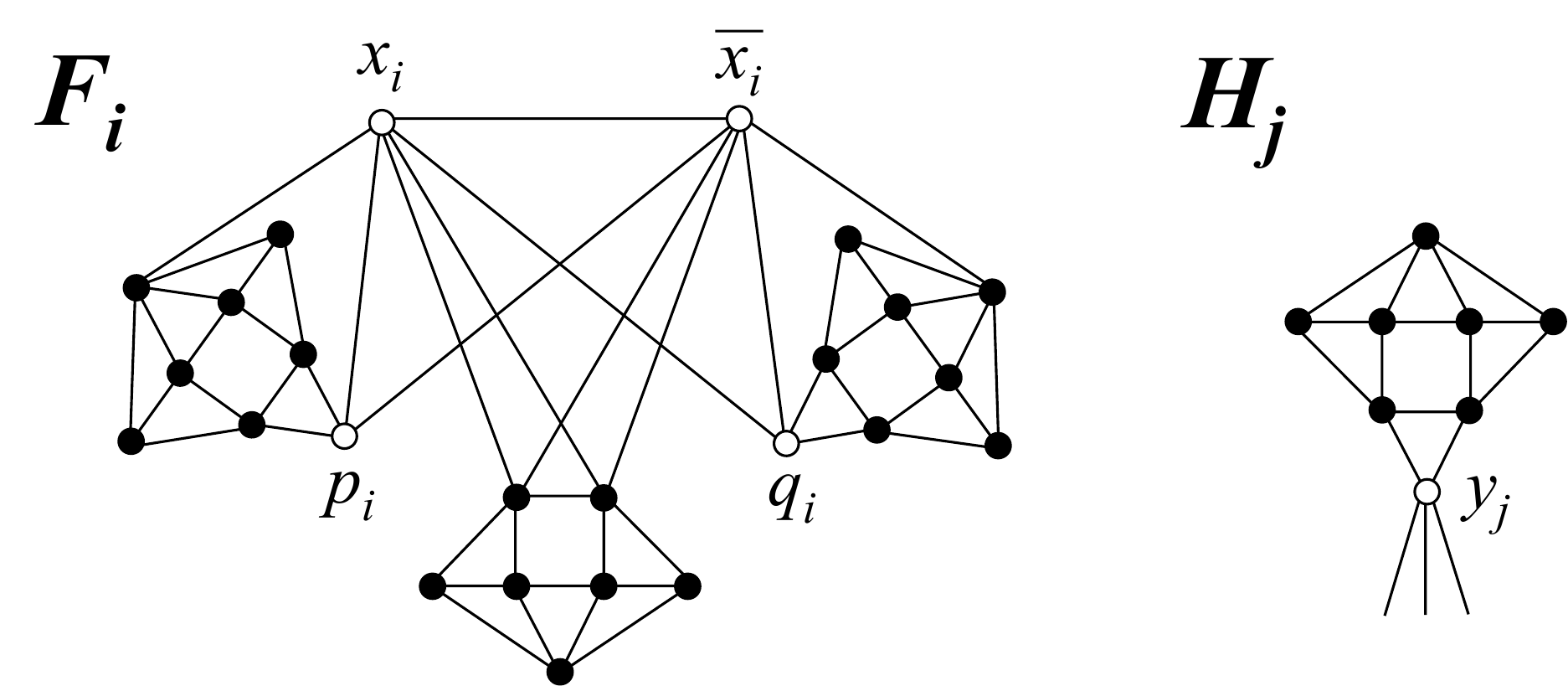}
    \caption{Variable and Clause graphs}
    \label{fig:variable-clause}
\end{wrapfigure}
\cbeginproof
\ref{lem:g7} 
Clearly, ERR:OLD is NP, as every possible candidate solution can be generated nondeterministically in polynomial time (specifically, $O(n)$ time), and each candidate can be verified in polynomial time using Theorem~\ref{theo:multi-param}~(\toroman{4}).
To complete the proof, we will now show a reduction from 3-SAT to ERR:OLD.

\begin{figure}[ht]
    \centering
    \includegraphics[width=0.7\textwidth]{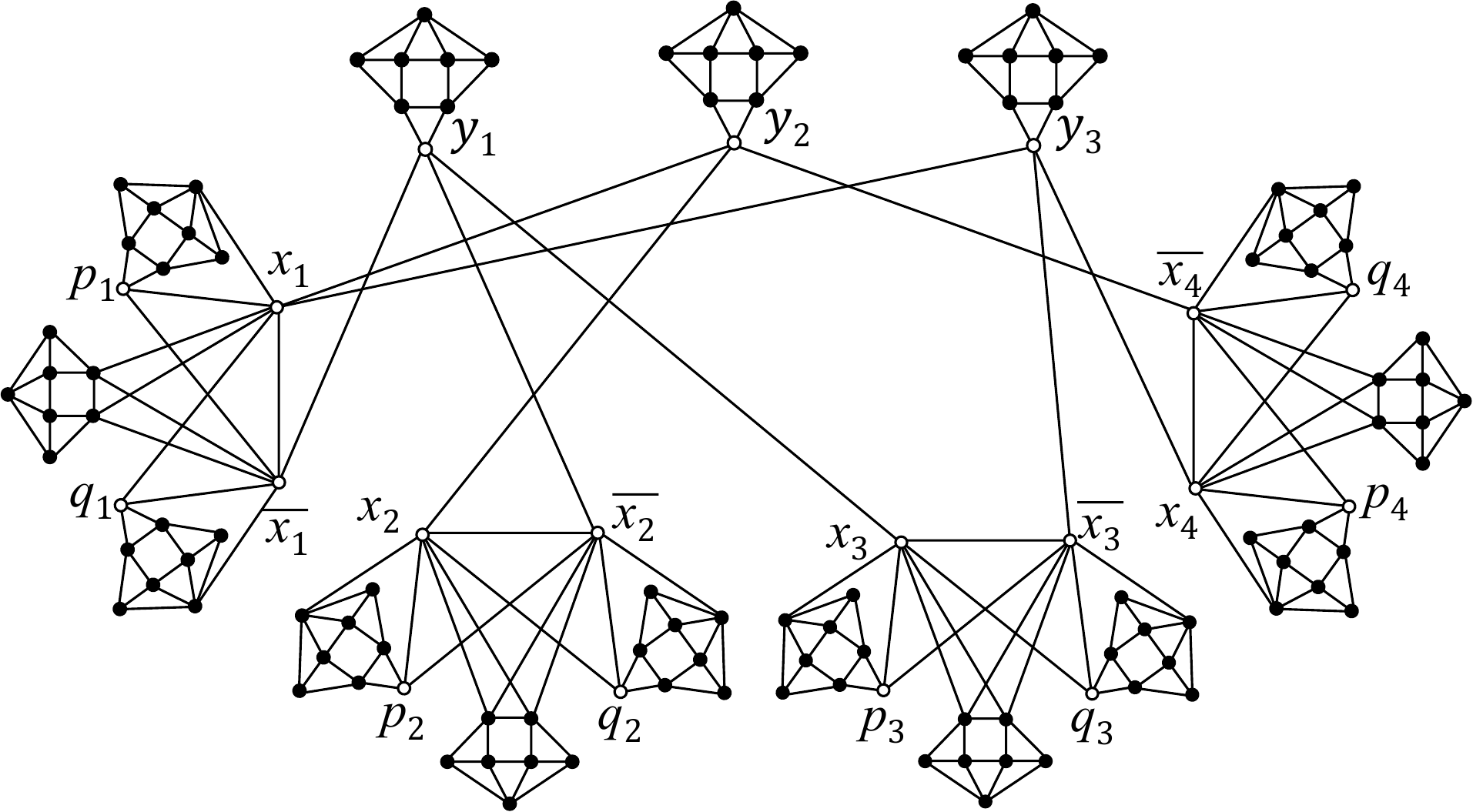}
    \caption{\begin{tabular}[t]{l} Construction of $G$ from $(\overline{x_1} \lor \overline{x_2} \lor x_3) \land (x_1 \lor x_2 \lor \overline{x_4}) \land (x_1 \lor \overline{x_3} \lor x_4)$ \protect\\ with $N = 4$, $M = 3$, $K = 109$  \end{tabular}}
    \label{fig:example-clause}
\end{figure}

Let $\psi$ be an instance of the 3-SAT problem with $M$ clauses on $N$ variables.
We will construct a graph, $G$, as follows.
For each variable $x_i$, create an instance of the $F_i$ graph (Figure~\ref{fig:variable-clause}); this includes a vertex for $x_i$ and its negation $\overline{x_i}$.
For each clause $c_j$ of $\psi$, create a new instance of the $H_j$ graph (Figure~\ref{fig:variable-clause}).
For each clause $c_j = \alpha \lor \beta \lor \gamma$, create an edge from the $c_j$ vertex to $\alpha$, $\beta$, and $\gamma$ from the variable graphs, each of which is either some $x_i$ or $\overline{x_i}$; for an example, see Figure~\ref{fig:example-clause}.
The resulting graph has precisely $25N + 8M$ vertices and $51N + 17M$ edges, and can be constructed in polynomial time.
To complete the problem instance, we define $K = 22N + 7M$.

Suppose $S \subseteq V(G)$ is an ERR:OLD set on $G$ with $|S| \le K$.
By Lemma~\ref{lem:g7}, we require at least all of the $21N + 7M$ detectors shown by the shaded vertices in Figure~\ref{fig:variable-clause}.
Additionally, in each $F_i$ the vertices $p_i$ and $q_i$ require $x_i \in S$ or $\overline{x_i} \in S$ in order to be 3-dominated.
Thus, we find that $|S| \ge 22N + 7M = K$, implying that $|S| = K$, so $|\{x_i,\overline{x_i}\} \cap S| = 1$ for each $i \in \{1, \hdots, N \}$.
For each $H_j$, we see that $y_j$ is not 3-dominated unless it is adjacent to at least one additional detector from its three neighbors in the $F_i$ graphs; therefore, $\psi$ is satisfiable.

For the converse, suppose we have a solution to the 3-SAT problem $\psi$; we will show that there is a ERR:OLD, $S$, on $G$ with $|S| \le K$.
We construct $S$ by first including all of the $21N + 7M$ vertices as shown in Figure~\ref{fig:variable-clause}.
Next, for each variable, $x_i$, if $x_i$ is true then we let the vertex $x_i \in S$; otherwise, we let $\overline{x_i} \in S$.
Thus, the fully-constructed $S$ has $|S| = 22N + 7M = K$.
It can be shown that every vertex pair is now 3-distinguished, and that every vertex which is not a $y_j$ clause vertex is 3-dominated.
Because we selected each $x_i \in S$ or $\overline{x_i} \in S$ based on a satisfying truth assignment for $\psi$, each $y_j$ must be adjacent to at least one additional detector vertex from the $F_i$ graphs.
Thus, all vertices are 3-dominated, implying $S$ is an ERR:OLD set for $G$ with $|S| \le K$, completing the proof.
\cendproof

\section{Infinite Grids}\label{sec:inf-grid}

\begin{figure}[ht]
    \centering
    \begin{tabular}{c@{\hskip 2em}c}
        \includegraphics[width=0.42\textwidth]{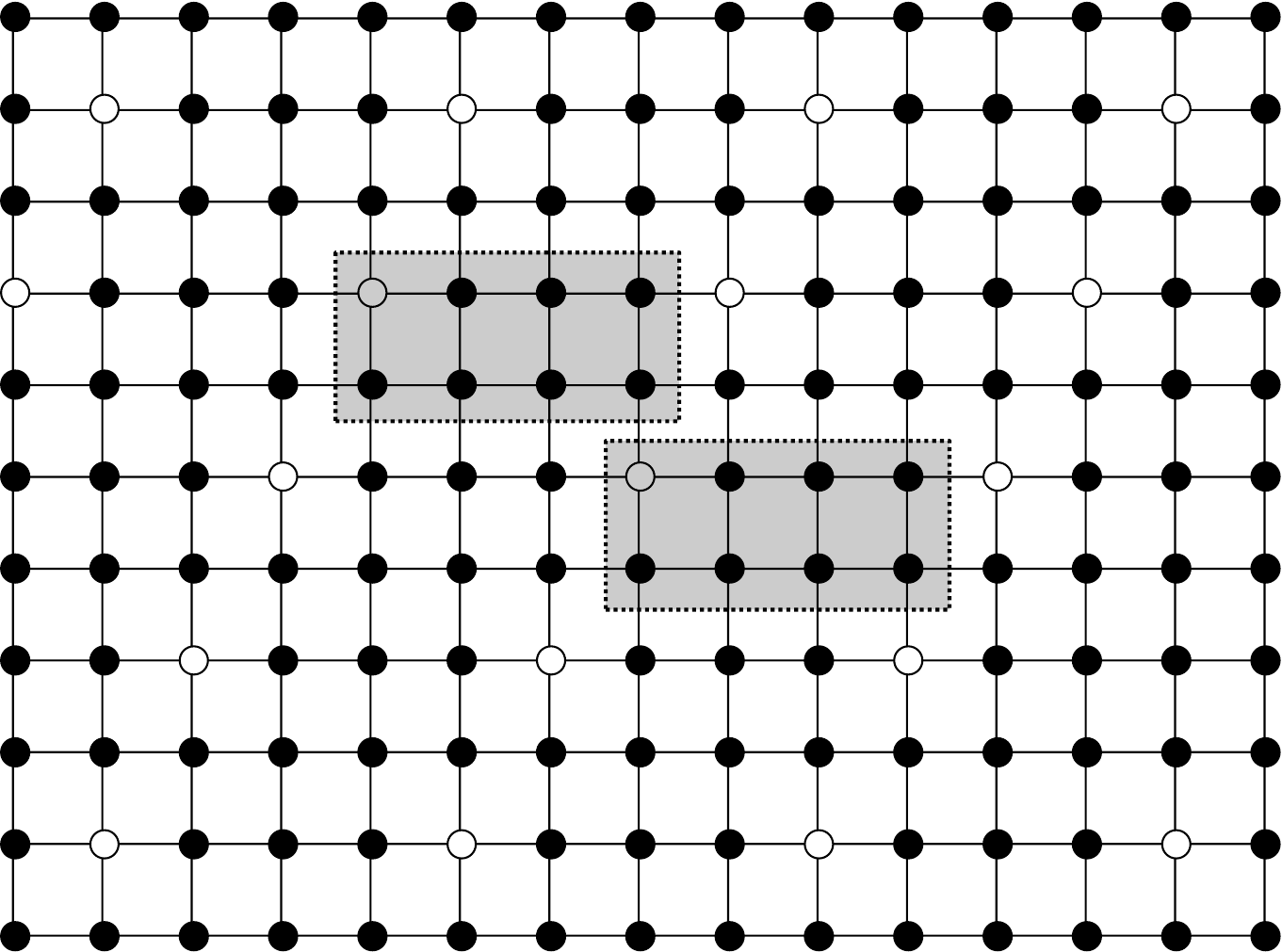} & \includegraphics[width=0.45\textwidth]{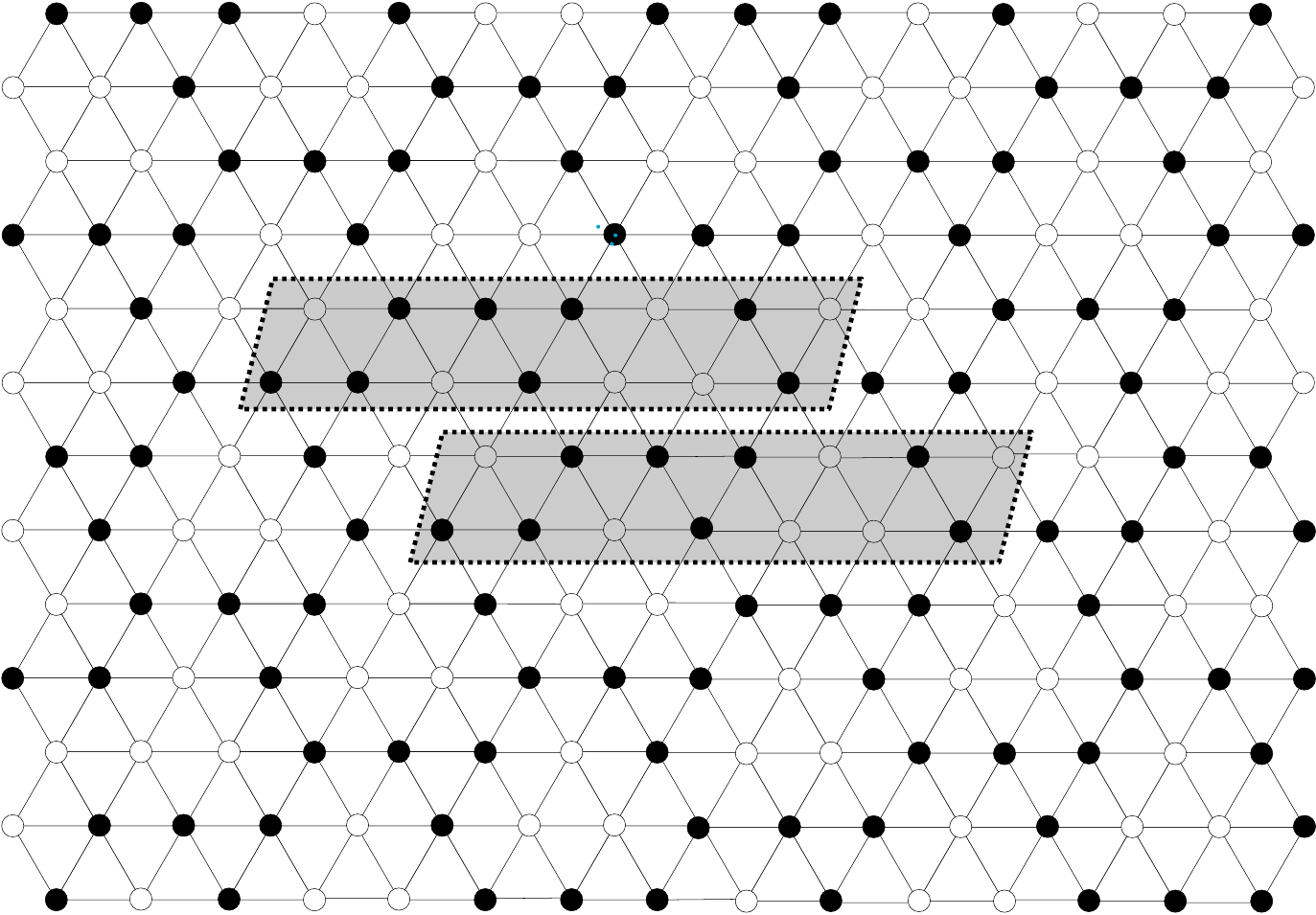} \\ (a) & (b) \\ \\
      \end{tabular}   \includegraphics[width=0.45\textwidth]{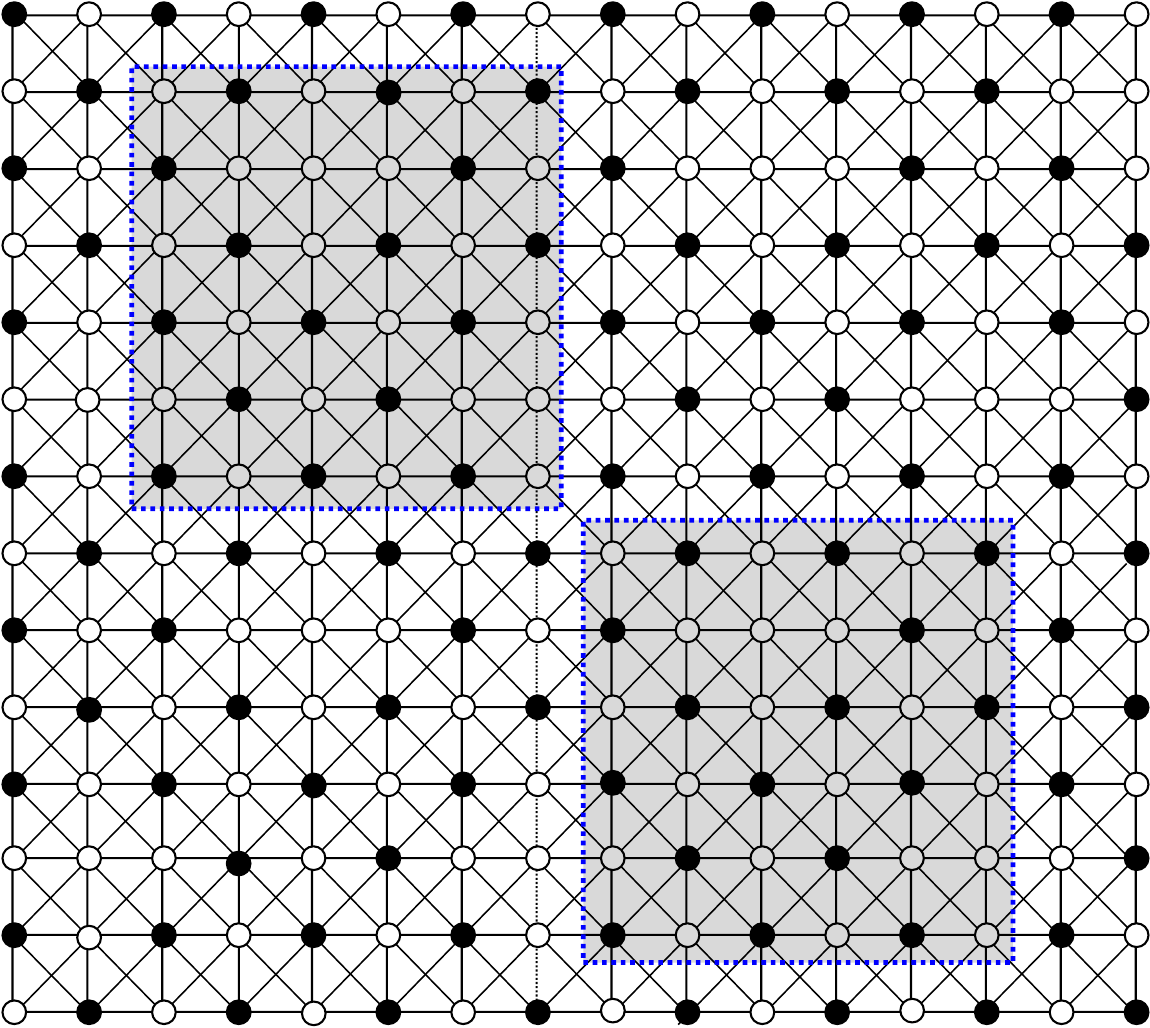}  \\ (c) 
   
    \caption{Our best constructions of ERR:OLD sets on SQR (a), TRI (b), and KNG (c). Shaded vertices denote detectors.}
    \label{fig:inf-grids-det-old-solns}
\end{figure}

The following theorem gives bounds for minimum ERR:OLD set densities in several infinite grids; the solutions achieving each upper bound can be found in Figure~\ref{fig:inf-grids-det-old-solns}.
Note that we do not include the long proofs for the lower bounds in Theorem~\ref{theo:grid-bounds}, but they can be demonstrated using a conventional share argument with discharging \cite{our-17}.

\begin{theorem}\label{theo:grid-bounds}
The upper and lower bounds on ERR:OLD:
\begin{enumerate}[label=\roman*]
    \item For the infinite square grid SQR, $\frac{6}{7} \le \textrm{ERR:OLD\%}(SQR) \le \frac{7}{8}$.
    \item For the infinite triangular grid TRI, $\frac{6}{11} \le \textrm{ERR:OLD\%}(TRI) \le \frac{4}{7}$.
    \item For the infinite king grid KNG,
$\frac{36}{83} \le \textrm{ERR:OLD\%}(KNG) \le \frac{4}{9}$.
\end{enumerate}
\end{theorem}

\FloatBarrier

\bibliographystyle{ACM-Reference-Format}
\bibliography{refs, lob-refs, extra-refs}


\begin{thebibliography}{14}


\ifx \showCODEN    \undefined \def \showCODEN     #1{\unskip}     \fi
\ifx \showDOI      \undefined \def \showDOI       #1{#1}\fi
\ifx \showISBNx    \undefined \def \showISBNx     #1{\unskip}     \fi
\ifx \showISBNxiii \undefined \def \showISBNxiii  #1{\unskip}     \fi
\ifx \showISSN     \undefined \def \showISSN      #1{\unskip}     \fi
\ifx \showLCCN     \undefined \def \showLCCN      #1{\unskip}     \fi
\ifx \shownote     \undefined \def \shownote      #1{#1}          \fi
\ifx \showarticletitle \undefined \def \showarticletitle #1{#1}   \fi
\ifx \showURL      \undefined \def \showURL       {\relax}        \fi
\providecommand\bibfield[2]{#2}
\providecommand\bibinfo[2]{#2}
\providecommand\natexlab[1]{#1}
\providecommand\showeprint[2][]{arXiv:#2}

\bibitem[Charon et~al\mbox{.}(2003)]%
        {ld-ic-np-complete-2}
\bibfield{author}{\bibinfo{person}{Irene Charon}, \bibinfo{person}{Olivier
  Hudry}, {and} \bibinfo{person}{Antoine Lobstein}.}
  \bibinfo{year}{2003}\natexlab{}.
\newblock \showarticletitle{{Minimizing the size of an identifying or
  locating-dominating code in a graph is NP-hard}}.
\newblock \bibinfo{journal}{\emph{Theoretical Computer Science}}
  \bibinfo{volume}{290}, \bibinfo{number}{3} (\bibinfo{year}{2003}),
  \bibinfo{pages}{2109--2120}.
\newblock
\urldef\tempurl%
\url{https://doi.org/10.1016/S0304-3975(02)00536-4}
\showDOI{\tempurl}


\bibitem[Cohen et~al\mbox{.}(2001)]%
        {NP-complete-ic}
\bibfield{author}{\bibinfo{person}{G. Cohen}, \bibinfo{person}{I. Honkala},
  \bibinfo{person}{A. Lobstein}, {and} \bibinfo{person}{G. Zémor}.}
  \bibinfo{year}{2001}\natexlab{}.
\newblock \showarticletitle{{On identifying codes}}, In
  \bibinfo{booktitle}{{Proc. DIMACS Workshop on Codes and Association
  Schemes}}.
\newblock \bibinfo{journal}{\emph{Discrete Mathematics and Theoretical Computer
  Science}}  \bibinfo{volume}{56}, \bibinfo{pages}{97--109}.
\newblock


\bibitem[Colbourn et~al\mbox{.}(1987)]%
        {NP-complete-ld}
\bibfield{author}{\bibinfo{person}{Charles~J Colbourn},
  \bibinfo{person}{Peter~J Slater}, {and} \bibinfo{person}{Lorna~K Stewart}.}
  \bibinfo{year}{1987}\natexlab{}.
\newblock \showarticletitle{{Locating dominating sets in series parallel
  networks}}.
\newblock \bibinfo{journal}{\emph{Congr. Numer}} \bibinfo{volume}{56},
  \bibinfo{number}{1987} (\bibinfo{year}{1987}), \bibinfo{pages}{135--162}.
\newblock


\bibitem[Dohner and Seo(2022)]%
        {redold}
\bibfield{author}{\bibinfo{person}{Robert Dohner} {and}
  \bibinfo{person}{Suk~Jai Seo}.} \bibinfo{year}{2022}\natexlab{}.
\newblock \showarticletitle{{The NP-completeness of Redundant
  Open-Locating-Dominating Set}}.
\newblock \bibinfo{journal}{\emph{CoRR}}  \bibinfo{volume}{abs/2201.05252}
  (\bibinfo{year}{2022}), \bibinfo{numpages}{9}~pages.
\newblock
\urldef\tempurl%
\url{https://doi.org/10.48550/arXiv.2201.05252}
\showDOI{\tempurl}
\showeprint[arXiv]{2201.05252}


\bibitem[Garey and Johnson(1979)]%
        {np-complete-bible}
\bibfield{author}{\bibinfo{person}{Michael~R. Garey} {and}
  \bibinfo{person}{David~S. Johnson}.} \bibinfo{year}{1979}\natexlab{}.
\newblock \bibinfo{booktitle}{\emph{{Computers and intractability: A guide to
  the theory of NP-completeness}}}.
\newblock \bibinfo{publisher}{W.H. Freeman}, \bibinfo{address}{San Francisco}.
\newblock
\showISBNx{978-0-7167-1045-5}


\bibitem[Honkala et~al\mbox{.}(2002)]%
        {honk02d}
\bibfield{author}{\bibinfo{person}{I. Honkala}, \bibinfo{person}{T. Laihonen},
  {and} \bibinfo{person}{S. Ranto}.} \bibinfo{year}{2002}\natexlab{}.
\newblock \showarticletitle{On strongly identifying codes}.
\newblock \bibinfo{journal}{\emph{Discrete Mathematics}}  \bibinfo{volume}{254}
  (\bibinfo{year}{2002}), \bibinfo{pages}{191--205}.
\newblock
\urldef\tempurl%
\url{https://doi.org/10.1016/S0012-365X(01)00357-0}
\showDOI{\tempurl}


\bibitem[Jean and Seo(2022)]%
        {our-17}
\bibfield{author}{\bibinfo{person}{Devin Jean} {and} \bibinfo{person}{Suk
  Seo}.} \bibinfo{year}{2022}\natexlab{}.
\newblock \bibinfo{title}{{Error-correcting Identifying Codes}}.
\newblock
\newblock
\showeprint[arxiv]{2204.11362}~[cs.DM]


\bibitem[Jean and Seo(2023a)]%
        {our-18}
\bibfield{author}{\bibinfo{person}{Devin Jean} {and} \bibinfo{person}{Suk
  Seo}.} \bibinfo{year}{2023}\natexlab{a}.
\newblock \bibinfo{title}{On Error-detecting Open-locating-dominating sets}.
\newblock
\newblock
\showeprint[arxiv]{2306.12583}~[math.CO]


\bibitem[Jean and Seo(2023b)]%
        {ourtri}
\bibfield{author}{\bibinfo{person}{Devin Jean} {and} \bibinfo{person}{Suk
  Seo}.} \bibinfo{year}{2023}\natexlab{b}.
\newblock \showarticletitle{{Optimal Error-Detecting Open-Locating-Dominating
  Set on the Infinite Triangular Grid}}.
\newblock \bibinfo{journal}{\emph{Discussiones Mathematicae: Graph Theory}}
  \bibinfo{volume}{43}, \bibinfo{number}{2} (\bibinfo{year}{2023}),
  \bibinfo{pages}{445--455}.
\newblock
\urldef\tempurl%
\url{https://doi.org/10.7151/dmgt.2374}
\showDOI{\tempurl}


\bibitem[Karpovsky et~al\mbox{.}(1998)]%
        {karpovsky}
\bibfield{author}{\bibinfo{person}{Mark~G Karpovsky},
  \bibinfo{person}{Krishnendu Chakrabarty}, {and} \bibinfo{person}{Lev~B
  Levitin}.} \bibinfo{year}{1998}\natexlab{}.
\newblock \showarticletitle{{On a new class of codes for identifying vertices
  in graphs}}.
\newblock \bibinfo{journal}{\emph{IEEE transactions on information theory}}
  \bibinfo{volume}{44}, \bibinfo{number}{2} (\bibinfo{year}{1998}),
  \bibinfo{pages}{599--611}.
\newblock
\urldef\tempurl%
\url{https://doi.org/10.1109/18.661507}
\showDOI{\tempurl}


\bibitem[lobstein(2022)]%
        {dombib}
\bibfield{author}{\bibinfo{person}{A. lobstein}.}
  \bibinfo{year}{2022}\natexlab{}.
\newblock \bibinfo{title}{Watching Systems, Identifying, Locating-Dominating
  and Discriminating Codes in Graphs}.
\newblock
  \bibinfo{howpublished}{\url{https://dragazo.github.io/bibdom/main.pdf}}.
\newblock


\bibitem[Seo and Slater(2015)]%
        {ftold}
\bibfield{author}{\bibinfo{person}{Suk Seo} {and} \bibinfo{person}{Peter
  Slater}.} \bibinfo{year}{2015}\natexlab{}.
\newblock \showarticletitle{Fault Tolerant Detectors for Distinguishing Sets in
  Graphs}.
\newblock \bibinfo{journal}{\emph{Discussiones Mathematicae Graph Theory}}
  \bibinfo{volume}{35} (\bibinfo{year}{2015}), \bibinfo{pages}{797--818}.
\newblock
\urldef\tempurl%
\url{https://doi.org/10.7151/dmgt.1838}
\showDOI{\tempurl}


\bibitem[Seo and Slater(2010)]%
        {old}
\bibfield{author}{\bibinfo{person}{Suk~J Seo} {and} \bibinfo{person}{Peter~J
  Slater}.} \bibinfo{year}{2010}\natexlab{}.
\newblock \showarticletitle{{Open neighborhood locating dominating sets}}.
\newblock
  \bibinfo{howpublished}{\url{https://ajc.maths.uq.edu.au/pdf/46/ajc_v46_p109.pdf}}.
\newblock \bibinfo{journal}{\emph{Australas. J Comb.}}  \bibinfo{volume}{46}
  (\bibinfo{year}{2010}), \bibinfo{pages}{109--120}.
\newblock


\bibitem[Slater(1987)]%
        {dom-loc-acyclic}
\bibfield{author}{\bibinfo{person}{Peter~J Slater}.}
  \bibinfo{year}{1987}\natexlab{}.
\newblock \showarticletitle{{Domination and location in acyclic graphs}}.
\newblock \bibinfo{journal}{\emph{Networks}} \bibinfo{volume}{17},
  \bibinfo{number}{1} (\bibinfo{year}{1987}), \bibinfo{pages}{55--64}.
\newblock


\end{thebibliography}

\end{document}